# On binary relations without non-identical endomorphisms

Apoloniusz Tyszka

**Summary**. On every set $A$ there is a rigid binary relation, i.e. such a relation $R$ that there is no homomorphism $(A,R) \to (A,R)$ except the identity (Vopěnka et al. [1965]). We state two conjectures which strengthen this theorem.

CONJECTURE 1. If $\kappa$ is an infinite cardinal number and card $A \leqslant 2^{2^\kappa}$ then there exists a relation $R \subseteq A \times A$ which satisfies

$$\forall_{\substack{x,y \in A \\ x \neq y}} \exists_{\substack{\{x\} \subseteq A(x,y) \subseteq A \\ \text{card } A(x,y) \leqslant \kappa}} \forall_{\substack{f: A(x,y) \to A \\ f(x)=y}} f \text{ is not a homomorphism of } R.$$

CONJECTURE 2. If $\kappa \neq 0$ is a limit cardinal number and card $A \leqslant 2^{\sup\{2^\alpha : \alpha \in \text{Card}, \alpha < \kappa\}}$ then there exists a relation $R \subseteq A \times A$ which satisfies

$$\forall_{\substack{x,y \in A \\ x \neq y}} \exists_{\substack{\{x\} \subseteq A(x,y) \subseteq A \\ \text{card } A(x,y) < \kappa}} \forall_{\substack{f: A(x,y) \to A \\ f(x)=y}} f \text{ is not a homomorphism of } R.$$

Conjecture 2 is valid for $\kappa = \omega$ (Tyszka [1994]). We prove Conjecture 1 for $\kappa = \omega$.

On every set $A$ there is a rigid binary relation, i.e. such a relation $R$ that there is no homomorphism $(A,R) \to (A,R)$ except the identity (see [4] and [1]). The following two conjectures were stated in [2] in a more complicated form.

---

**Conjecture 1**. If $\kappa$ is an infinite cardinal number and card $A \leq 2^{2^\kappa}$ then there exists a relation $R \subseteq A \times A$ which satisfies the following condition (∗):

$$(*) \quad \forall \begin{array}{c} x,y \in A \\ x \neq y \end{array} \exists \begin{array}{c} \{x\} \subseteq A(x,y) \subseteq A \\ \text{card } A(x,y) \leq \kappa \end{array} \forall \begin{array}{c} f: A(x,y) \to A \\ f(x)=y \end{array} f \text{ is not a homomorphism of } R.$$

**Remark 1** ([2]). If $R \subseteq A \times A$ satisfies the condition (∗) then $R$ is rigid. If $\kappa$ is an infinite cardinal number and a relation $R \subseteq A \times A$ satisfies the condition (∗) then card $A \leq 2^{2^\kappa}$.

**Conjecture 2**. If $\kappa \neq 0$ is a limit cardinal number and card $A \leq 2^{\sup\{2^\alpha : \alpha \in \text{Card}, \alpha < \kappa\}}$ then there exists a relation $R \subseteq A \times A$ which satisfies the following condition (∗∗):

$$(**) \quad \forall \begin{array}{c} x,y \in A \\ x \neq y \end{array} \exists \begin{array}{c} \{x\} \subseteq A(x,y) \subseteq A \\ \text{card } A(x,y) < \kappa \end{array} \forall \begin{array}{c} f: A(x,y) \to A \\ f(x)=y \end{array} f \text{ is not a homomorphism of } R.$$

**Remark 2** ([2]). If $R \subseteq A \times A$ satisfies the condition (∗∗) then $R$ is rigid. If $\kappa \neq 0$ is a limit cardinal number and a relation $R \subseteq A \times A$ satisfies the condition (∗∗) then card $A \leq 2^{\sup\{2^\alpha : \alpha \in \text{Card}, \alpha < \kappa\}}$.

**Theorem 1** ([2]). Conjecture 2 is valid for $\kappa = \omega$.

In this paper we prove Conjecture 1 for $\kappa = \omega$. We begin from introductory results. Let $p_0 = (0,0), p_1 = (1,0), p_2 = (\frac{1}{2}, \frac{\sqrt{3}}{2})$, $U := \{(a,b) \in \mathbb{R}^2 : |a-b| = 1\}$, $\Phi := \{R \subseteq \mathbb{R}^2 \times \mathbb{R}^2 : R \cup R^{-1} = U, R \cap R^{-1} = \{(p_0, p_1), (p_1, p_0), (p_0, p_2), (p_2, p_0)\}, (p_1, p_2) \in R\}$.

**Observation 1**. If $S, Z \in \Phi$, $\{p_0, p_1, p_2\} \subseteq X \subseteq \mathbb{R}^2$, $f: X \to \mathbb{R}^2$ and $f: \langle X, S \rangle \to \langle \mathbb{R}^2, Z \rangle$ is a homomorphism, then $f(p_0) = p_0$, $f(p_1) = p_1$, $f(p_2) = p_2$.

**Observation 2**. If $S, Z \in \Phi$, $X \subseteq \mathbb{R}^2$, $f: X \to \mathbb{R}^2$ and $f: \langle X, S \rangle \to \langle \mathbb{R}^2, Z \rangle$ is a homomorphism, then $f$ is a homomorphism of $S \cup S^{-1} = Z \cup Z^{-1} = U$ i.e. $f$ preserves all unit distances.

**Lemma** ([3]). If $x, y \in \mathbb{R}^n$ ($n > 1$) and $\varepsilon > 0$ then there exists a finite set $T_{xy}(\varepsilon) \subseteq \mathbb{R}^n$ containing $x$ and $y$ such that each map $f: T_{xy}(\varepsilon) \to \mathbb{R}^n$ preserving all unit distances satisfies $||f(x) - f(y)| - |x-y|| \leq \varepsilon$.

**Note**. Let us consider the finite set $T_{xy}(\varepsilon) \subseteq \mathbb{R}^2$ from the above Lemma. For each $t \in T_{xy}(\varepsilon)$ there exist points $t(0), t(1), \ldots, t(m(t)) \in \mathbb{R}^2$ ($m(t) \in \omega$) satisfying: $t(0)=t$, $t(m(t))=x$, $|t(i)-t(i+1)|=1$ ($0 \leq i \leq m(t)-1$). The finite set

$$\tilde{T}_{xy}(\varepsilon) := \bigcup_{t \in T_{xy}(\varepsilon)} \{t(0), t(1), \ldots, t(m(t))\}$$

has the property from the Lemma and the graph $U$ is connected on $\tilde{T}_{xy}(\varepsilon)$.

The following Theorem 2 was proved in [2] in a more complicated way. Theorem 2 immediately implies a weaker version of Theorem 1.

**Theorem 2**. If $R \in \Phi$ then $R$ satisfies the condition (**) for $\kappa = \omega$.

*Proof.* Assume that $R \in \Phi$, $x, y \in \mathbb{R}^2$ and $x \neq y$. Since $p_0, p_1, p_2$ are not collinear there exists $i \in \{0, 1, 2\}$ satisfying $|p_i - x| \neq |p_i - y|$. Let

$$\mathbb{R}^2(x, y) := \{p_0, p_1, p_2\} \cup T_{p_i x}\left(\frac{||p_i - x| - |p_i - y||}{2}\right).$$

The set $\mathbb{R}^2(x, y)$ is finite. Assume that $f: \mathbb{R}^2(x, y) \to \mathbb{R}^2$ is a homomorphism of $R$. We need to prove that $f(x) \neq y$. According to Observation 1 $f(p_i) = p_i$, by Observation 2 $f$ preserves all unit distances. From this, applying the Lemma we obtain

$$|p_i - f(x)| = |f(p_i) - f(x)| \in \left[|p_i - x| - \frac{||p_i - x| - |p_i - y||}{2}, \ |p_i - x| + \frac{||p_i - x| - |p_i - y||}{2}\right].$$

On the other hand

$$|p_i - y| \notin \left[|p_i - x| - \frac{||p_i - x| - |p_i - y||}{2}, \ |p_i - x| + \frac{||p_i - x| - |p_i - y||}{2}\right].$$

Therefore $f(x) \neq y$. This completes the proof.

We are now in a position to formulate the main result.

**Theorem 3**. If $\emptyset \neq \Psi \subseteq \Phi$ then the relation $R_\Psi \subseteq (\mathbb{R}^2 \times \Psi) \times (\mathbb{R}^2 \times \Psi)$ defined by the following formula

$$\forall \ x, y \in \mathbb{R}^2 \ \forall \ S, Z \in \Psi \ \left(((x, S), (y, Z)) \in R_\Psi \Leftrightarrow (x, y) \in S = Z\right)$$

satisfies the condition (*) for $\kappa = \omega$.

*Proof.* Assume that $(x,S), (y,Z) \in \mathbb{R}^2 \times \Psi$ and $(x,S) \neq (y,Z)$. Let $\pi: \{S\} \to \{Z\}$.

First case: $x \neq y$. Since $p_0, p_1, p_2$ are not collinear there exists $i \in \{0,1,2\}$ satisfying $|p_i - x| \neq |p_i - y|$. Let

$$(\mathbb{R}^2 \times \Psi)((x,S),(y,Z)) := \left(\{p_0, p_1, p_2\} \cup \tilde{T}_{p_i x}\left(\frac{||p_i - x| - |p_i - y||}{2}\right)\right) \times \{S\}.$$

The set $(\mathbb{R}^2 \times \Psi)((x,S),(y,Z))$ is finite. Suppose, on the contrary, that $f: (\mathbb{R}^2 \times \Psi)((x,S),(y,Z)) \to \mathbb{R}^2 \times \Psi$ is a homomorphism of $R_\Psi$ and $f((x,S)) = (y,Z)$. By the Note the graph $R_\Psi \cup R_\Psi^{-1}$ is connected on $(\mathbb{R}^2 \times \Psi)((x,S),(y,Z))$. Therefore $f$ maps

$$\left(\{p_0, p_1, p_2\} \cup \tilde{T}_{p_i x}\left(\frac{||p_i - x| - |p_i - y||}{2}\right)\right) \times \{S\}$$

into $\mathbb{R}^2 \times \{Z\}$. From this, there is a uniquely determined transformation

$$\bar{f} : \{p_0, p_1, p_2\} \cup \tilde{T}_{p_i x}\left(\frac{||p_i - x| - |p_i - y||}{2}\right) \to \mathbb{R}^2$$

satisfying: $f = \langle \bar{f}, \pi \rangle$, $\bar{f}(x) = y$. Obviously

$$\bar{f} : \left\langle \{p_0, p_1, p_2\} \cup \tilde{T}_{p_i x}\left(\frac{||p_i - x| - |p_i - y||}{2}\right), S \right\rangle \to \left\langle \mathbb{R}^2, Z \right\rangle$$

is a homomorphism. According to Observation 1 $\bar{f}(p_i) = p_i$, by Observation 2 $\bar{f}$ preserves all unit distances. From this, applying the Note and the Lemma we obtain

$$|p_i - \bar{f}(x)| = |\bar{f}(p_i) - \bar{f}(x)| \in \left[|p_i - x| - \frac{||p_i - x| - |p_i - y||}{2}, |p_i - x| + \frac{||p_i - x| - |p_i - y||}{2}\right].$$

On the other hand

$$|p_i - y| \notin \left[|p_i - x| - \frac{||p_i - x| - |p_i - y||}{2}, |p_i - x| + \frac{||p_i - x| - |p_i - y||}{2}\right].$$

Therefore $\bar{f}(x) \neq y$. This contradicts our assumption.

Second case: $x = y$ so $S \neq Z$. Since $S \neq Z$ there exist $u, v \in \mathbb{R}^2$ such that $|u - v| = 1$, $(u,v) \in S$, $(v,u) \notin S$, $(v,u) \in Z$, $(u,v) \notin Z$. There exist points $x_0, x_1, \ldots, x_k \in \mathbb{R}^2$ ($k \in \omega$) satisfying: $x_0 = x$, $x_k = p_0$, $|x_i - x_{i+1}| = 1$ ($0 \leq i \leq k-1$). Let

$$(\mathbb{R}^2 \times \Psi)((x,S),(y,Z)) := \left(\{x_0, x_1, \ldots, x_k\} \cup \bigcup_{\substack{j \in \omega \\ i \in \{0,1,2\}}} \tilde{T}_{p_i u}\left(\frac{1}{j}\right) \cup \tilde{T}_{p_i v}\left(\frac{1}{j}\right)\right) \times \{S\}.$$

- 4 -

The set $(\mathbb{R}^2 \times \Psi)((x,S),(y,Z))$ is countable. Suppose, on the contrary, that $f:(\mathbb{R}^2 \times \Psi)((x,S),(y,Z)) \to \mathbb{R}^2 \times \Psi$ is a homomorphism of $R_\Psi$ and $f((x,S))=(y,Z)$. By the Note the graph $R_\Psi \cup R_\Psi^{-1}$ is connected on $(\mathbb{R}^2 \times \Psi)((x,S),(y,Z))$. Therefore $f$ maps

$$\left( \{x_0, x_1, \ldots, x_k\} \cup \bigcup_{\substack{j \in \omega \\ i \in \{0,1,2\}}} \tilde{T}_{p_i u}\left(\frac{1}{j}\right) \cup \tilde{T}_{p_i v}\left(\frac{1}{j}\right) \right) \times \{S\}$$

into $\mathbb{R}^2 \times \{Z\}$. From this, there is a uniquely determined transformation

$$\bar{f} : \{x_0, x_1, \ldots, x_k\} \cup \bigcup_{\substack{j \in \omega \\ i \in \{0,1,2\}}} \tilde{T}_{p_i u}\left(\frac{1}{j}\right) \cup \tilde{T}_{p_i v}\left(\frac{1}{j}\right) \to \mathbb{R}^2$$

satisfying: $f=\langle \bar{f}, \pi \rangle$, $\bar{f}(x)=y$. Obviously

$$\bar{f} : \left\langle \{x_0, x_1, \ldots, x_k\} \cup \bigcup_{\substack{j \in \omega \\ i \in \{0,1,2\}}} \tilde{T}_{p_i u}\left(\frac{1}{j}\right) \cup \tilde{T}_{p_i v}\left(\frac{1}{j}\right) , S \right\rangle \to \langle \mathbb{R}^2, Z \rangle$$

is a homomorphism. By Observation 2 $\bar{f}$ preserves all unit distances. Therefore, by the Note and the Lemma $\bar{f}$ preserves distances $|p_0-u|$, $|p_1-u|$, $|p_2-u|$, $|p_0-v|$, $|p_1-v|$, $|p_2-v|$. On the other hand, according to Observation 1 $\bar{f}(p_0)=p_0$, $\bar{f}(p_1)=p_1$, $\bar{f}(p_2)=p_2$. From this $\bar{f}(u)=u$ and $\bar{f}(v)=v$. Since $f$ is a homomorphism of $R_\Psi$ and $(u,v) \in S$ we conclude that $((\bar{f}(u), \pi(S)), (\bar{f}(v), \pi(S))) \in R_\Psi$ i.e. $(u,v) \in Z$. This contradicts our assumption.

We have proved Theorem 3.

**Corollary**. Conjecture 1 is valid for $\kappa = \omega$.

*Proof*. Assume that $\kappa = \omega$. According to Theorem 3 if $2^\omega \leq \operatorname{card} A \leq 2^{2^\omega}$ then there exists a relation $R \subseteq A \times A$ which satisfies the condition (∗). On the other hand, in virtue of Theorem 1 if $\operatorname{card} A \leq 2^{\sup\{2^\alpha : \alpha \in \operatorname{Card}, \alpha < \omega\}} = 2^\omega$ then there exists a relation $R \subseteq A \times A$ which satisfies the stronger condition (∗∗), so the condition (∗) holds too.

A. Tyszka
Technical Faculty
Hugo Kollataj University
Balicka 104
PL-30-149 Krakow
Poland
E-mail: *rttyszka@cyf-kr.edu.pl*